\documentclass[hidelinks,12pt]{amsart}
\usepackage{extsizes}
\usepackage{bbm}
\usepackage{libertine}
\usepackage{enumitem}
\usepackage[T1]{fontenc}
\usepackage{stmaryrd}
\usepackage{color}
\DeclareFontFamily{OT1}{pzc}{}
\DeclareFontShape{OT1}{pzc}{m}{it}{<-> s * [1.100] pzcmi7t}{}
\DeclareMathAlphabet{\mathpzc}{OT1}{pzc}{m}{it}\DeclareMathAlphabet{\mathcal}{OMS}{cmsy}{m}{n}\usepackage{fullpage}
\usepackage{mdframed}
\usepackage{graphicx}
\usepackage{calligra}
\usepackage{wasysym}
\usepackage{csquotes}
\usepackage{bm}
\usepackage{amsmath}
\usepackage[all]{xy}
\usepackage{amssymb}
\usepackage{amsthm}
\usepackage{mathrsfs}
\usepackage[OT2,T1]{fontenc}
\usepackage{centernot}
\usepackage{hyperref}

\newcommand{\DMO}[2]{\DeclareMathOperator{#1}{#2}}\DMO{\BK}{BK}
\DMO{\FL}{FL}
\DMO{\Ann}{Ann}
\DMO{\std}{std}
\DMO{\antidiag}{antidiag}
\DMO{\locadm}{loc.adm}
\DMO{\Inj}{Inj}
\DMO{\LL}{LL}
\DMO{\Dmod}{\emph{D }-mod}
\DMO{\univ}{univ}
\DMO{\Fitt}{Fitt}
\DMO{\WD}{WD}
\DMO{\geom}{geom}
\DMO{\Fl}{Fl}
\DMO{\grad}{grad}
\DMO{\labmda}{\lambda}
\DMO{\Iw}{Iw}
\DMO{\tor}{tor}
\DMO{\coh}{coh}
\DMO{\vol}{vol}
\DMO{\semsim}{ss}
\DMO{\free}{free}
\DMO{\Alg}{Alg}
\DMO{\oth}{otherwise}
\DMO{\Ber}{Ber}
\DMO{\Diff}{Diff}
\DMO{\br}{br}
\DMO{\Isot}{Isot}
\DMO{\prim}{prim}
\DMO{\RAH}{RAH}
\DMO{\Sets}{Sets}
\DMO{\cone}{cone}
\DMO{\Grps}{Grps}
\DMO{\Dec}{Dec}
\DMO{\Flat}{Flat}
\DMO{\AbGps}{AbGps}
\DMO{\Sch}{Sch}
\DMO{\AH}{AH}
\DMO{\cl}{cl}
\DMO{\sk}{sk}
\DMO{\HC}{HC}
\DMO{\cosk}{sk}
\DMO{\ur}{ur}
\DMO{\LocSys}{LocSys}
\DMO{\rk}{rk}
\DMO{\NT}{NT}
\DMO{\cork}{cork}
\DMO{\KS}{KS}
\DMO{\MU}{MU}
\DMO{\der}{der}
\DMO{\Art}{Art}
\DMO{\Proj}{Proj}
\DMO{\End}{End}
\DMO{\Betti}{Betti}
\DMO{\Sym}{Sym}
\DMO{\cInd}{cInd}
\DMO{\GL}{GL}
\DMO{\Gal}{Gal}
\DMO{\Br}{Br}
\DMO{\Der}{Der}
\DMO{\Sp}{Sp}
\DMO{\Tan}{Tan}
\DMO{\Spin}{Spin}
\DMO{\Var}{Var}
\DMO{\Nrd}{Nrd}
\DMO{\cusp}{cusp}
\DMO{\Mat}{Mat}
\DMO{\Isom}{Isom}
\DMO{\Stab}{Stab}
\DMO{\SO}{SO}
\DMO{\Res}{Res}
\DMO{\Lie}{Lie}
\DMO{\SU}{SU}
\DMO{\Ad}{Ad}
\DMO{\ad}{ad}
\DMO{\im}{im}
\DMO{\Frob}{Frob}
\DMO{\Fr}{Fr}
\DMO{\red}{red}
\DMO{\an}{an}
\DMO{\Pic}{Pic}
\DMO{\Tor}{Tor}
\DMO{\Hdg}{Hdg}
\DMO{\id}{id}
\DMO{\pr}{pr}
\DMO{\Mor}{Mor}
\DMO{\Ext}{Ext}
\DMO{\ML}{ML}
\DMO{\PGL}{PGL}
\DMO{\SL}{SL}
\DMO{\GU}{GU}
\DMO{\GSp}{GSp}
\DMO{\GSL}{GSL}
\DMO{\Aff}{Aff}
\DMO{\NS}{NS}
\DMO{\gr}{gr}
\DMO{\Ch}{Ch}
\DMO{\QCoh}{QCoh}
\DMO{\Coh}{Coh}
\DMO{\inv}{inv}
\DMO{\Gr}{Gr}
\DMO{\Bun}{Bun}
\DMO{\Hk}{Hk}
\DMO{\GH}{GH}
\DMO{\HT}{HT}
\DMO{\LT}{LT}
\DMO{\Int}{Int}
\DMO{\UU}{U}
\DMO{\OO}{O}
\DMO{\Loc}{Loc}
\DMO{\Conn}{Conn}
\DMO{\sing}{sing}
\DMO{\si}{si}
\DMO{\Sen}{Sen}
\DMO{\MaxSpec}{MaxSpec}
\DMO{\ran}{ran}
\DMO{\coker}{coker}
\DMO{\DIV}{div}
\DMO{\Cl}{Cl}
\DMO{\Frac}{Frac}
\DMO{\VEC}{Vec}
\DMO{\Weil}{Weil}
\DMO{\SPLIT}{split}
\DMO{\Tr}{Tr}
\DMO{\val}{val}
\DMO{\pv}{p.v.}
\DMO{\disc}{disc}
\DMO{\trdeg}{tr.deg}
\DMO{\rad}{rad}
\DMO{\codim}{codim}
\DMO{\dist}{dist}
\DMO{\length}{length}
\DMO{\diam}{diam}
\DMO{\Supp}{Supp}
\DMO{\Ass}{Ass}
\DMO{\ord}{ord}
\DMO{\RE}{Re}
\DMO{\Sh}{Sh}
\DMO{\IM}{Im}
\DMO{\Tot}{Tot}
\DMO{\Bl}{Bl}
\DMO{\lcm}{lcm}
\DMO{\ann}{ann}
\DMO{\arcsinh}{arcsinh}
\DMO{\CHAR}{char}
\DMO{\MOD}{mod}
\DMO{\BB}{BB}
\DMO{\new}{new}
\DMO{\alg}{alg}
\DMO{\Irr}{Irr}
\DMO{\res}{res}
\DMO{\rank}{rank}
\DMO{\naive}{naive}
\DMO{\tors}{tors}
\DMO{\Perf}{Perf}
\DMO{\Sht}{Sht}
\DMO{\Perv}{Perv}
\DMO{\soc}{soc}
\DMO{\Mod}{Mod}
\DMO{\cyc}{cyc}
\DMO{\SC}{sc}
\DMO{\SP}{sp}
\DMO{\Deck}{Deck}
\DMO{\PSL}{PSL}
\DMO{\Area}{Area}
\DMO{\Cont}{Cont}
\DMO{\sgn}{sgn}
\DMO{\Cat}{Cat}
\DMO{\Cov}{Cov}
\DMO{\rig}{rig}
\DMO{\FSch}{FSch}
\DMO{\Rig}{Rig}
\DMO{\Spv}{Spv}
\DMO{\Spa}{Spa}
\DMO{\trace}{trace}
\DMO{\cont}{cont}
\DMO{\aff}{aff}
\DMO{\cor}{cor}
\DMO{\CH}{CH}
\DMO{\Spec}{Spec}
\DMO{\rec}{rec}
\DMO{\LGC}{LGC}
\DMO{\un}{un}
\DMO{\conj}{conj}
\DMO{\Eval}{Eval}
\DMO{\JH}{JH}
\DMO{\can}{can}
\DMO{\Fss}{Fss}
\DMO{\Speh}{Speh}
\DMO{\Ind}{Ind}
\DMO{\ch}{ch}
\DMO{\nr}{nr}
\DMO{\Swan}{Swan}
\DMO{\St}{St}
\DMO{\Ho}{Ho}
\DMO{\HH}{HH}
\DMO{\trop}{trop}
\DMO{\Jac}{Jac}
\DMO{\vir}{vir}
\DMO{\coll}{coll}
\DMO{\reg}{reg}
\DMO{\dlog}{dlog}
\DMO{\Div}{Div}
\DMO{\ab}{ab}
\DMO{\Tam}{Tam}
\DMO{\Ran}{Ran}
\DMO{\IC}{IC}
\DMO{\Sat}{Sat}
\DMO{\Rat}{Rat}
\DMO{\loc}{loc}
\DMO{\ev}{ev}
\DMO{\st}{st}
\DMO{\pst}{pst}
\DMO{\Fil}{Fil}
\DMO{\cris}{cris}
\DMO{\dR}{dR}
\DMO{\Rep}{Rep}
\DMO{\Sel}{Sel}
\DMO{\spec}{spec}
\DMO{\Spf}{Spf}
\DMO{\JL}{JL}
\DMO{\BGL}{BGL}
\DMO{\Arc}{Arc}
\DMO{\MHS}{MHS}
\DMO{\Nm}{Nm}
\DMO{\holim}{holim}
\DMO{\nInd}{nInd}
\DMO{\sSets}{s\textbf{Sets}}
\DMO{\sArt}{s\textbf{Art}}
\DMO{\BDJ}{BDJ}
\DMO{\GV}{GV}
\DMO{\BM}{BM}
\DMO{\Ord}{Ord}
\DMO{\mult}{mult}
\DMO{\WDRep}{WDRep}
\DMO{\Aut}{Aut}
\DMO{\Hom}{Hom}
\DMO{\sph}{sph}
\DMO{\Def}{Def}
\DMO{\GO}{GO}
\DMO{\diag}{diag}
\DMO{\cond}{cond}
\DMO{\ind}{ind}
\DMO{\irr}{irr}
\DMO{\RHom}{RHom}
\DMO{\sm}{sm}
\DMO{\sss}{ss}
\DMO{\sHom}{sHom}
\DMO{\Tran}{Tran}
\DMO{\Rees}{Rees}
\DMO{\lcv}{lcv} 
\DMO{\SN}{SN}
\DMO{\triv}{triv}
\DMO{\height}{ht}
\DMO{\proj}{proj}
\DMO{\Fun}{Fun}
\DMO{\cts}{cts}
\DMO{\Obj}{Obj}
\DMO{\Sing}{Sing}
\DMO{\Pro}{Pro}
\DMO{\Ig}{Ig}
\DMO{\Ha}{Ha}
\DMO{\BC}{BC}
\DMO{\RZ}{RZ}
\DMO{\supp}{supp}
\DMO{\projdim}{proj.dim}
\DMO{\Zar}{Zar}
\DMO{\Ban}{Ban}
\DMO{\LA}{LA}
\DMO{\ess}{ess}
\DMO{\op}{op}
\DMO{\Func}{Func}
\DMO{\Born}{Born}
\DMO{\Comm}{Comm}
\DMO{\Dr}{Dr}
\DMO{\LC}{LC}
\DMO{\nind}{n-ind}
\DMO{\perf}{perf}
\DMO{\charpoly}{char.poly}
\makeatletter
\def\thmhead@plain#1#2#3{%
  \thmname{#1}\thmnumber{\@ifnotempty{#1}{ }\@upn{#2}}%
  \thmnote{ {\the\thm@notefont#3}}}
\let\thmhead\thmhead@plain
\makeatother

\newtheorem*{thm1*}{Theorem}
\newtheorem*{lemma*}{Lemma}
\newtheorem*{defn1*}{Definition}

\newtheorem{thm2}{Theorem}[section]

\newtheorem*{prop*}{Proposition}

\newtheorem{prop2}[thm2]{Proposition}

\newtheorem*{conj*}{Conjecture}

\theoremstyle{definition}

\newtheorem*{defn2*}{Definition}

\newtheorem{homework}{}
\newtheorem*{prb*}{Problem}

\newtheorem*{claim*}{Claim}

\newtheorem*{rmk*}{Remark}

\newtheorem{exam2}[thm2]{Example}

\newtheoremstyle{theoremdd}
  {6pt}
  {6pt}
  {}
  {0pt}
  {\bfseries}
  {.}
  { }
  {\thmname{#1}\thmnumber{ #2}\textnormal{\thmnote{ #3}}}
  \theoremstyle{theoremdd}

\newtheoremstyle{theoremee}
  {6pt}
  {6pt}
  {\itshape}
  {0pt}
  {\bfseries}
  {.}
  { }
  {\thmname{#1}\thmnumber{ #2}\textnormal{\thmnote{ #3}}}
  \theoremstyle{theoremee}

\newcommand{\xrar}[1]{\xrightarrow{#1}}

\newcommand{\riso}{\xrar{\sim}}

 \newenvironment{psmat}
  {\left(\begin{smallmatrix}}
  {\end{smallmatrix}\right)}

 \newenvironment{psmatrix}
  {\left(\begin{smallmatrix}}
  {\end{smallmatrix}\right)}

\newcommand{\textbox}[2] {\left\lbrace\parbox{#1em}{\center{#2}}\right\rbrace}

\newcommand{\wt}{\widetilde}

\newcommand{\ov}{\overline}

\newcommand{\rar}{\rightarrow}

\newcommand{\ncom}[1]{\newcommand{#1}}
\ncom{\sbuset}{\subset}

\newcommand{\hrar}{\hookrightarrow}

\newcommand{\und}[1]{\underline{#1}}
\newcommand{\emphC}[1]{\textsf{\textbf{#1}}}

\DeclareSymbolFont{cyrletters}{OT2}{wncyr}{m}{n}
\DeclareMathSymbol{\Sha}{\mathalpha}{cyrletters}{"58}

\newcommand{\et}{\operatorname{\acute{e}t}}

\DMO{\bmr}{\mathbbm{r}}
\DMO{\bmf}{\mathbbm{f}}
\DMO{\bmx}{\mathbbm{x}}

\newcommand{\bC}{\mathbb{C}}

\newcommand{\bF}{\mathbb{F}}

\newcommand{\bH}{\mathbb{H}}

\newcommand{\bQ}{\mathbb{Q}}

\newcommand{\bZ}{\mathbb{Z}}

\newcommand{\sT}{\mathscr{T}}

\newcommand{\cL}{\mathcal{L}}

\newcommand{\cO}{\mathcal{O}}
\newcommand{\cP}{\mathcal{P}}

\newcommand{\cT}{\mathcal{T}}

\newcommand{\fE}{\mathfrak{E}}

\newcommand{\fS}{\mathfrak{S}}
\newcommand{\fT}{\mathfrak{T}}
\newcommand{\fU}{\mathfrak{U}}
\newcommand{\fV}{\mathfrak{V}}

\newcommand{\MSRI}{\let\thefootnote\relax\footnotetext{A part of the material is based upon work supported by the National Science Foundation under Grant No. DMS-1928930 while the author was in residence at the Mathematical Sciences Research Institute in Berkeley, California during the Spring 2023 semester.}}

\DMO{\Vect}{Vect}
\DMO{\Out}{Out}
\DMO{\sub}{sub}
\DMO{\quo}{quo}
\begin{document}
\title{A proof of N\'eron--Ogg--Shafarevich criterion via its archimedean analogue}
\author{Gyujin Oh}
\maketitle
\begin{abstract}
In this short note, we deduce the classical N\'eron--Ogg--Shafarevich criterion on good reduction of abelian varieties from its \emph{archimedean analogue}: a holomorphic family of abelian varieties over a punctured disc extends to the whole unit disc if and only if the topological monodromy representation is trivial.
\end{abstract}
\tableofcontents
\section{Introduction}
There is a folklore analogy between the ramification theory of $\ell$-adic Galois representations and the singularity of vector bundles with integrable connections. A general dictionary is that the unramified, tamely ramified, and wildly ramified Galois representations correspond to the removable, regular, and irregular singularities, respectively. Furthermore, an integrable connection with at worst regular singularities is encoded by its topological monodromy representation by the Riemann--Hilbert corresopndence. Thus, one expects that the unramifiedness of $\ell$-adic monodromy representation is related to the triviality of topological monodromy representation.

On the other hand, the classical N\'eron--Ogg--Shafarevich criterion relates the good reduction of an abelian variety with the unramifiedness of its $\ell$-adic monodromy representation. Guided by the viewpoint of the previous paragraph, we first formulate and prove the corresponding criterion regarding the triviality of topological monodromy representation of a family of abelian varieties, which we call the \emph{archimedean analogue} of the N\'eron--Ogg--Shafarevich criterion. To avoid confusion, from now on, we will refer to the classical N\'eron--Ogg--Shafarevich criterion as the $\ell$-adic N\'eron--Ogg--Shafarevich criterion.

The main purpose of the note is to deduce the $\ell$-adic N\'eron--Ogg--Shafarevich criterion from the archimedean N\'eron--Ogg--Shafarevich criterion. The proof crucially uses the arithmetic toroidal compactification of the moduli of principally polarized abelian varieties as in the work of Faltings--Chai \cite{FC}. This new proof of the $\ell$-adic criterion explicitly realizes the folklore analogy between the ramification theory of $\ell$-adic Galois representations and the singularities of topological monodromy representations. 

The idea of the proof is to use the ``degeneration of degenerations.'' Namely, we use the fact that Mumford's degeneration of abelian varieties can realize all semi-abelian degeneration of abelian varieties, and also that such degeneration can be put into a family. Perhaps the most well-known instance of this fact is that, for elliptic curves, the Tate uniformization can be put into a universal family, called the \emph{universal Tate elliptic curve} of Raynaud. The idea of relating $\ell$-adic and topological monodromy via a family of degenerations was prominently used in \cite{Oda}, which studies the Galois representation of a higher-genus curve using the information about its topological monodromy. 

The article is organized as follows. In \S2, we will first precisely formulate and prove the archimedean N\'eron--Ogg--Shafarevich criterion. Then, in \S3, we will first deduce the $\ell$-adic N\'eron--Ogg--Shafarevich criterion from the archimedean criterion in the case of elliptic curves, using the more familiar universal Tate elliptic curve. In \S4, we will end with the proof of the $\ell$-adic criterion in the general case by using Mumford's construction of degeneration of abelian varieties. 
\subsection*{Acknowledgements} We thank John Halliday for pointing out a mistake in a previous version.
\section{Archimedean N\'eron--Ogg--Shafarevich criterion}
We first formulate the ``archimedean analogue'' of the N\'eron--Ogg--Shafarevich criterion, which is considerably easier to prove than the $\ell$-adic N\'eron--Ogg--Shafarevich criterion.
\begin{prop2}[(Archimedean N\'eron--Ogg--Shafarevich criterion)]\label{ArchNOS}Let $f:A\rar D^{\times}$ be a holomorphic family of abelian varieties of dimension $g$ over the punctured disc $D^{\times}$. Then, $f$ extends to a family of abelian varieties over $D$ if and only if the monodromy representation $\rho:\pi_{1}(D^{\times},t_{0})\cong\bZ\rar\Aut H^{1}(A_{t_{0}},\bZ)$ is trivial, where $t_{0}\in D^{\times}$ is a fixed base point.
\end{prop2}
\begin{proof}
If the family extends to the unit disc $D$, the monodromy representation factors through $\pi_{1}(D,t_{0})$, which is trivial. Conversely, suppose the monodromy representation is trivial. We first prove this direction assuming that the family $f:A\rar D^{\times}$ is a family of \emph{principally polarized} abelian varieties. Then, the family $f$ defines a period morphism from $D^{\times}$ modulo the monodromy to the Siegel upper half plane $\bH_{g}$. As the monodromy is trivial, the period morphism is a holomorphic map $p:D^{\times}\rar\bH_{g}$. As $\bH_{g}$ is conformally equivalent to a bounded domain, $0=D-D^{\times}$ is a removable singularity of $p$.

Now let us prove the converse direction in the general case. By Zahrin's trick, there is a family of principally polarized abelian varieties $A'\rar D^{\times}$ whose topological monodromy representation is trivial and which contains $A\subset A'$ as a family of abelian subvarieties. By the previous paragraph, $A'\rar D^{\times}$ extends to a family of principally polarized abelian varieties $\wt{A'}\rar D$. Let $\wt{A}$ be the closure of $A$ in $\wt{A'}$, with its reduced scheme structure. We claim that $\wt{A}\rar D$ is a family of abelian varieties, or that $\wt{A}_{0}$ is an abelian variety. Note that the fiberwise group structure $\mu:\wt{A'}\times_{D}\wt{A'}\rar \wt{A'}$ restricts to $\mu':\wt{A}\times_{D}\wt{A}\rar\wt{A'}$. Since $\mu'\rvert_{A\times_{D^{\times}}A}$ factors through $A\subset\wt{A'}$, $\mu'$ factors through $\wt{A}$. This shows that $\wt{A}_{0}$ is a reduced closed algebraic subgroup of the abelian variety $\wt{A'}_{0}$. Since a reduced connected closed algebraic subgroup of an abelian variety over a characteristic zero field is an abelian variety, it is now sufficient to show that $\wt{A}_{0}$ is connected, which is true as it is the continuous image of the closure of a connected set $A\times_{D^{\times}}A$, which is itself connected.
\end{proof}
Using this, we would like to prove the $\ell$-adic N\'eron--Ogg--Shafarevich criterion. From now on, we fix a finite extension $K/\bQ_{p}$.
\begin{thm2}[($\ell$-adic N\'eron--Ogg--Shafarevich criterion)]\label{NOS}Let $A/K$ be an abelian variety. Let $\rho:G_{K}\rar\GL(T_{\ell}A)$ be the $\ell$-adic monodromy representation for $\ell\ne p$. Then, $A$ has good reduction if and only if $\rho\rvert_{I_{K}}$ is trivial.
\end{thm2}
The rest of the article will be dedicated to proving Theorem \ref{NOS}. As one direction is immediate, we are left to prove the other nontrivial direction, namely proving the good reduction of $A$ assuming that $\rho\rvert_{I_{K}}$ is trivial. 

As the first reduction step, we note that both sides of the statement of Theorem \ref{NOS} are invariant under an unramified base change. Furthermore, by Zahrin's trick, one can assume that $A$ is principally polarized. 

Next, we note that only semistable abelian varieties need to be considered, using the geometry of the arithmetic toroidal compactification of the moduli of abelian varieties with level structures.
\begin{prop2}\label{Reduction1}
Let $A/K$ be a principally polarized abelian variety, and $\rho$ be the $\ell$-adic monodromy representation of $A$ for $\ell\ne p$. If $\rho\rvert_{I_{K}}$ is trivial, then $A$ has at worst semistable reduction.
\end{prop2}
\begin{proof}
Let $g$ be the dimension of $A$. By assumption, $A[\ell^{n}]$ is unramified for any $n\ge1$. Take a finite unramified extension $L/K$ over which $A[\ell^{3}]$ splits. Then, $A[\ell^{3}](L)\cong(\bZ/\ell^{3}\bZ)^{2g}$. Fix one such isomorphism $\iota$. Then, $(A,\iota)$ defines a point $P\in A_{g,\ell^{3}}(L)$, where $A_{g,\ell^{3}}$ is the moduli space of principally polarized abelian varieties of dimension $g$ with full level $\ell^{3}$ structure. By \cite[Theorem 7.10]{GIT}, it is known that $A_{g,\ell^{3}}$ is represented by a smooth quasi-projective scheme over $\bZ[1/\ell]$. Furthermore, by \cite[Theorem IV.6.7]{FC}, for a good enough choice of auxiliary data, there is a smooth projective $\bZ[1/\ell]$-scheme $\ov{A}_{g,n}$, which contains $A_{g,n}$ as a dense open subscheme. By the valuative criterion for properness, the point $P\in A_{g,\ell^{3}}(L)$ extends to an $\cO_{L}$-point $\cP\in\ov{A}_{g,\ell^{3}}(\cO_{L})$.

Recall that there is a semi-abelian scheme $\ov{G}\rar\ov{A}_{g,\ell^{3}}$, extending the universal abelian scheme $G\rar A_{g,n}$. After pulling back this family via the $\cO_{L}$-point $\cP:\Spec\cO_{L}\rar \ov{A}_{g,\ell^{3}}$, we obtain a semi-abelian scheme over $\cO_{L}$ whose generic fiber is $A_{L}$. This implies that $A_{L}$ has at worst semistable reduction. As the reduction type is invariant under an unramified base change, we conclude that $A$ has at worst semistable reduction, as desired.
\end{proof}\section{Warm-up: the case of elliptic curves}

Mumford's construction of degenerating abelian varieties in the case of elliptic curves is exactly the \emph{Tate uniformization}. By the means of Raynaud's universal Tate elliptic curve, one even knows that  the degeneration of elliptic curves can be put into a one-dimensional family. In this section, we show how to prove Theorem \ref{NOS} using the universal Tate elliptic curve, to illustrate the usefulness of a family of degenerations.
\begin{proof}[Proof of Theorem \ref{NOS}, in the case of elliptic curves]
Suppose that an elliptic curve $A$ over $K$ has bad reduction but also has trivial monodromy over the inertia group $I_{K}$. As we know from Proposition \ref{Reduction1}, $A$ has semistable reduction. We will leverage the fact that $A$ has the Tate uniformization, and furthermore that the Tate uniformization is obtained from the \emph{universal Tate elliptic curve} of Raynaud. Our reference for the construction is \cite[\S9.2]{Bosch} and \cite[\S2.5]{Conrad}.

Let $K^{\nr}$ be the maximal unramified extension of $K$. Then, the \emph{universal Tate elliptic curve} is a family $f:\fT\rar\fS=\Spec\cO_{K^{\nr}}[[q]]$ 
 such that, if $E$ is an elliptic curve over $K^{\nr}$ with semistable bad reduction, then the (semistable) N\'eron model $\fE$ of $E$ over $\Spec \cO_{K^{\nr}}$ 
 is obtained by pulling back $\fT$ along the map $\cO_{K^{\nr}}[[q]]\rar\cO_{K^{\nr}}$, $q\mapsto q(E)$, the $q$-parameter of $E$. The family $\fT$ used here is obtained by taking the base-change of $\und{\mathrm{Tate}}_{1}$ in \cite[\S2.5]{Conrad}, which is in turn the algebraization of the (formal) universal Tate elliptic curve in \cite[\S9.2]{Bosch}, from $\bZ[[q]]$ to $\cO_{K^{\nr}}[[q]]$.

Let $\fV=\Spec\cO_{K^{\nr}}$ be the closed subscheme of $\fS$ defined by the 
map $\cO_{K^{\nr}}[[q]]\xrar{q\mapsto q(A)}\cO_{K^{\nr}}$, and let $\eta$ be the generic point of $\fV$. Let $i:\fV\hrar\fS$ be the closed embedding. Then, \[f\rvert_{\fV}:i^{*}\fT\rar\fV,\]is a semistable model of $A_{K^{\nr}}$. This in particular means that $f\rvert_{\fV_{\eta}}:(i^{*}\fT)_{\eta}\rar\fV_{\eta}$ is identified with $A_{K^{\nr}}\rar\Spec(K^{\nr})$. Upon choosing a geometric generic point $\ov{\eta}:\Spec(\ov{K})\rar\Spec(K^{\nr})\rar\fV$, the pro-$\ell$-part of the monodromy representation \[\rho_{\fV_{\eta},\ov{\eta}}:\pi_{1,\et}(\fV_{\eta},\ov{\eta})_{\ell}\rar\Aut(R^{1}(f\rvert_{{\ov{\eta}}})_{*}\und{{\bQ}_{\ell}}),\]is identified with the $\ell$-adic monodromy representation of $A$ restricted to the pro-$\ell$-part of the inertia group $(I_{K})_{\ell}$, \[\rho\rvert_{(I_{K})_{\ell}}:(I_{K})_{\ell}\rar\Aut(H^{1}(A_{\ov{K}},{\bQ}_{\ell})).\]
On the other hand, as $A_{K^{\nr}}\rar\Spec(K^{\nr})$ is put in a larger family, the above monodromy representation $\rho_{\fV_{\eta},\ov{\eta}}$ factors through the monodromy over a larger base $\fU=\lbrace q\ne0\rbrace=\Spec\cO_{K^{\nr}}((q))\subset\fV$:
\[\xymatrix{\pi_{1,\et}(\fV_{\eta},\ov{\eta})_{\ell}\ar[rd]\ar[rr]^-{\rho_{\fV_{\eta},\ov{\eta}}}&&\Aut(R^{1}(f\rvert_{\fV_{\ov{\eta}}})_{*}\und{\bQ_{\ell}})\\&\pi_{1,\et}(\fU,\ov{\eta})_{\ell}\ar[ru]_-{\rho_{\fU,\ov{\eta}}}}\]namely, the left arrow is the natural map $\pi_{1,\et}(\fV_{\eta},\ov{\eta})_{\ell}\rar\pi_{1,\et}(\fU,\ov{\eta})_{\ell}$ via $\cO_{K^{\nr}}((q))\xrar{q\mapsto q(A)}K^{\nr}$. Note that, by Abhyankar's lemma, all \'etale $\ell$-covers of both $\fU$ and $\fV_{\eta}$ are Kummer covers, so both $\pi_{1,\et}(\fV_{\eta},\ov{\eta})_{\ell}$ and $\pi_{1,\et}(\fU,\ov{\eta})_{\ell}$ are isomorphic to $\bZ_{\ell}$. After identifying both with $\bZ_{\ell}$, the natural map is multiplication by $\ell^{n}$ for some $n\ge0$.\footnote{This was pointed out to us by John Halliday.} By assumption, $\rho_{\fV_{\eta},\ov{\eta}}=\rho\rvert_{(I_{K})_{\ell}}$ is trivial, so we conclude that $\rho_{\fU,\ov{\eta}}$ has finite image.

On the other hand, $\fT$ is smooth and proper over $\fU$, so the specialization map of \'etale fundamental groups identifies the monodromy over $\fU$ at $\ov{\eta}$ is identified with the monodromy over $\fU$ at $\ov{\iota}$, where $\ov{\iota}=\Spec \ov{K^{\nr}((q))}$ is a geometric generic point of $\fU=\Spec\cO_{K^{\nr}}((q))$. Thus, the monodromy representation
\[\rho_{\fU,\ov{\iota}}:\pi_{1,\et}(\fU,\ov{\iota})_{\ell}\rar\Aut(R^{1}(f\rvert_{\ov{\iota}})_{*}\und{\bQ_{\ell}}),\]has finite image.
Now $\pi_{1,\et}(\cO_{K^{\nr}})=1$, so $\rho_{\fU,\ov{\iota}}$ will stay the same even if we base change $\cO_{K^{\nr}}$ to a characteristic zero algebraically closed field. We fix an abstract field embedding $j:K_{\nr}\hrar\bC$ and consider the base-change to $\bC$ using this embedding. Therefore, the monodromy of the $\bC$-base-change
\[\rho_{\fU_{\bC},\ov{\iota}_{\bC}}:\pi_{1,\et}(\fU_{\bC},\ov{\iota}_{\bC})_{\ell}=\Gal(\ov{\bC((q))}/\bC((q)))_{\ell}\rar\Aut(R^{1}(f\rvert_{\ov{\iota}_{\bC}})_{*}\und{\bQ_{\ell}}),\]has finite image.

The crucial point is that the family $\fT_{\bC[[q]]}\rar\Spec\bC[[q]]$, obtained by base-changing $f:\fT\rar\fS$ via $\cO_{K^{\nr}}\hrar K^{\nr}\xrar{j}\bC$, is the formal germ of a complex-analytic family $f':\cT\rar D$ of complex analytic varieties over the unit disk. One could see this by noticing that the power series defining the Tate uniformization has coefficients in $\bZ$ and can be verbatim used to define a complex-analytic family over the punctured unit disk. Another way to realize this is that the aforementioned family induces a map $p:\Spec\bC[[q]]\rar X(1)$, where $X(1)$ is the moduli of elliptic curves. As $X(1)$ is an algebraic space, one could take the open unit disk around the point $p\in X(1)^{\an}$, where $X(1)^{\an}$ is a Moishezon manifold.

Using the above observation, we see that $\rho_{\fU_{\bC},\ov{\iota}_{\bC}}$ is the $\ell$-completion of the topological monodromy $\rho_{f'}:\pi_{1}(D^{\times},t_{0})\rar\Aut H^{1}(A_{t_{0}},\bZ)$, where $t_{0}\in D^{\times}$. As the family $\cT\rar D$ is semistable, this local topological monodromy is unipotent. However, any unipotent matrix in $\GL_{2}(\bZ)$ is either trivial or has infinite order, so this concludes that the topological monodromy over $D^{\times}$ is trivial. By the archimedean N\'eron--Ogg--Shafarevich criterion, Proposition \ref{ArchNOS}, $\cT\rvert_{D^{\times}}\rar D^{\times}$ can be extended to a family of abelian varieties $\cT'\rar D$. As $\cT\rvert_{D^{\times}}\rar D^{\times}$ can be extended to the semistable family $\cT\rar D$, this contradicts the fact that there is at most one way of extending a family of elliptic curves over $D^{\times}$ to $D$.
\end{proof}

\section{A proof of N\'eron--Ogg--Shafarevich criterion}
The crucial part of the proof of the previous section is that you can put a semistable elliptic curve into a family of semi-abelian varieties over $\Spec\cO_{K^{\nr}}[[q]]$ whose base-change to $\bC[[q]]$ is the formal germ of a  complex-analytic family of semi-abelian varieties. With this in mind, we can now prove Theorem \ref{NOS} in the general case.
\begin{proof}[Proof of Theorem \ref{NOS}]
As in the proof of Proposition \ref{Reduction1}, there is a point $P\in A_{g,\ell^{3}}(K^{\nr})$ such that $G_{P}\cong A_{K^{\nr}}$, and there is an extension $i:\Spec\cO_{K^{\nr}}\rar \ov{A}_{g,\ell^{3}}$. Let $\ov{P}\in A_{g,\ell^{3}}(\ov{\bF}_{p})$ be the point in the special fiber of $i$. As $\ov{A}_{g,\ell^{3}}$ is smooth over $\bZ[1/\ell]$, the local ring of $\ov{A}_{g,\ell^{3}}$ at $\ov{P}$ is \[\cO_{\ov{A}_{g,\ell^{3}},\ov{P}}\cong W(\ov{\bF}_{p})[[X_{1},\cdots,X_{d}]], \]where $d=\frac{g(g+1)}{2}$. The morphism $i:\Spec\cO_{K^{\nr}}\rar\ov{A}_{g,\ell^{3}}$ then induces $\Spf\cO_{K^{\nr}}\rar\Spf\cO_{\ov{A}_{g,\ell^{3}},\ov{P}}$, or a map $f:W(\ov{\bF}_{p})[[X_{1},\cdots,X_{d}]]\rar\cO_{K^{\nr}}$. As its mod $p$ reduction is $\ov{P}$, if we choose a uniformizer $\pi\in\cO_{K^{\nr}}$, then $f(X_{i})$ is divisible by $\pi$. Thus, $f$ factors through
\[W(\ov{\bF}_{p})[[X_{1},\cdots,X_{d}]]\xrar{X_{i}\mapsto \frac{f(X_{i})}{\pi}X}\cO_{K^{\nr}}[[X]]\xrar{X\mapsto \pi}\cO_{K^{\nr}}.\]Let $p:\Spec\cO_{K^{\nr}}[[X]]\rar \ov{A}_{g,\ell^{3}}$ be the corresponding morphism, and consider $p^{*}\ov{G}\rar\Spec\cO_{K^{\nr}}[[X]]$. This is, as before, a one-dimensional family of semi-abelian schemes, which has the following properties.\begin{itemize}\item It specializes to $A_{K^{\nr}}$ at a certain point.
\item It is smooth away from the $X=0$ locus.
\item Its base-change to $\bC[[X]]$ is the formal germ of a complex-analytic semistable family of abelian varieties over a unit disc.
\end{itemize}
The proof then proceeds exactly as in the proof in \S3, where the only difference is that we have to use the fact that a non-identity unipotent matrix in $\GL_{n}(\bZ)$ has infinite order for any $n$.
\end{proof}

\bibliographystyle{alpha}

\begin{thebibliography}{OOOO}
\bibitem[Bos]{Bosch} Siegfried Bosch, Lectures on formal and rigid geometry. Lecture Notes in Mathematics, 2105. \emph{Springer, Cham,} 2014.
\bibitem[Con]{Conrad} Brian Conrad, Arithmetic moduli of generalized elliptic curves. \emph{J. Inst. Math. Jussieu} 6 (2007), no. 2, 209-278.
\bibitem[FC]{FC} Gerd Faltings, Ching-Li Chai, Degeneration of abelian varieties. With an appendix by David Mumford. Ergebnisse der Mathematik und ihrer Grenzgebiete (3), 22. \emph{Springer-Verlag, Berlin}, 1990.
\bibitem[Mum]{GIT} David Mumford, John Fogarty, Frances Kirwan, Geometric invariant theory. Third edition. Ergebnisse der Mathematik und ihrer Grenzgebiete (2), 34. \emph{Springer-Verlag, Berlin,} 1994. 
\bibitem[Oda]{Oda} Takayuki Oda, A note on ramification of the Galois representation on the fundamental group of an algebraic curve. II. \emph{J. Number Theory} 53 (1995), no. 2, 342–355.
\bibitem[SGA1]{SGA1} Alexander Grothendieck, Rev\^etements \'etales et groupe fondamental. Troisi\`eme \'edition, corrig\'ee S\'eminaire de G\'eom\'etrie Alg\'ebrique, 1960/61. \emph{Institut des Hautes \'Etudes Scientifiques, Paris}, 1963.
\end{thebibliography}

\end{document}